\documentclass[letterpaper, 11pt]{amsart}

\usepackage{mathptmx}
\usepackage{ulem}
\usepackage{graphics}
\usepackage{amstext}
\usepackage{amsmath}
\usepackage{graphicx}
\usepackage[latin1]{inputenc}
\usepackage{amsfonts}
\usepackage{amssymb}
\usepackage{graphics}
\usepackage[margin=1in]{geometry}
\usepackage{xcolor}
\usepackage{calligra}
\usepackage{pbsi}
\usepackage[T1]{fontenc}
\usepackage{latexsym}
\usepackage[all]{xy}

\newtheorem{theorem}{Theorem}
\newtheorem{lemma}[theorem]{Lemma}

\newtheorem{definition}[theorem]{Definition}

\newcommand{\Sym}{\mathrm{Sym}\,}
\newcommand{\Metz}{\mathrm{Metz}\,}

\title{Global product structure for a space of special matrices}

\author{Horacio~Leyva}
\email{hleyva@mat.uson.mx}

\author{Francisco~A.~Carrillo}
\email{carrillo@mat.uson.mx}

\author{Baltazar~Aguirre-Hern\'andez}
\email{bahe@xanum.uam.mx}

\author{Jes\'us F. Espinoza}
\email{jesus.espinoza@mat.uson.mx}

\address{Mathematics Department, Universidad de Sonora, Hermosillo, Sonora, M\'exico}


\begin{document}
\begin{abstract}
The importance of the Hurwitz Metzler matrices and the Hurwitz symmetric matrices can be appreciated in different applications: communication networks, biology and economics are some of them. In this paper, we use an approach of differential topology for studying such matrices. Our results are as follows: the space of the $n\times n$ Hurwitz symmetric matrices has a product manifold structure given by the space of the $(n-1)\times(n-1)$ Hurwitz symmetric matrices and the euclidean space. Additionally we study the space of Hurwitz Metzler matrices and these ideas let us do an analysis of robustness of Hurwitz Metzler matrices. In particular, we study the Insulin Model as application.
\end{abstract}

\maketitle

\section{Introduction}
The Hurwitz symmetric matrices are symmetric matrices whose eigenvalues have negative real part. On the other hand, a real matrix $A\in\mathcal{M}_{n\times n}$ is a Metzler matrix if its off diagonal elements are non-negative. 
The symmetric matrices are important in the analysis of stability with the Lyapunov approach (see \cite{Gantmacher:87}) while the importance of the Hurwitz symmetric matrices can be appreciated in \cite{Gantmacher:87}, \cite{Golub:89}, \cite{Horn:85} and \cite{Horn:91}. On the other hand, the Metzler matrices have great importance in applications (see for instance \cite{arcat-sontag}, \cite{johnson}, \cite{meyn}, \cite{siljak}, and \cite{stipanovic}). Now we use an approach of differential topology for studying the spaces of the $n\times n$ Hurwitz symmetric matrices and the $n\times n$ Hurwitz Metzler matrices, which are denoted by $\Sym\mathcal{H}_{n\times n}$ and $\Metz\mathcal{H}_{n\times n}$, respectively. Such approach proved to be useful in the study of the spaces of stable polynomials (see \cite{aguirre1} and \cite{aguirre2}). Our result use an idea presented in \cite{narendra}, \cite{redheffer1}, \cite{redheffer2} and \cite{shorten}: let $A_n$  be a $n\times n$ invertible matrix that is partitioned as $$A_n=\begin{pmatrix}
A_{n-1} & b_{n-1} \\ c_{n-1}^T & d_{n-1}
\end{pmatrix} \mbox{ and } D_{n-1} := A_{n-1}- \frac{b_{n-1}c_{n-1}^T}{d_{n-1}}$$
where $A_{n-1}\in\mathcal{M}_{(n-1)\times(n-1)}$, $b_{n-1},\ c_{n-1}\in\mathbb{R}^{n-1}$, $d_{n-1}\in\mathbb{R}^-$. Then, we have the following results of Narendra and Shorten \cite{narendra} and \cite{shorten}:

\begin{lemma}\label{lema1}
\begin{enumerate}
\item[1)] $A_n$ is a Hurwitz symmetric matrix if and only if $D_{n-1}$ is a Hurwitz symmetric matrix,
\item[2)] Let $A$ be a Metzler matrix. Then, $A$ is Hurwitz if and only if $D_{n-1}$ is a Hurwitz matrix.
\end{enumerate}
\end{lemma}

\noindent\textbf{Proof.} See \cite{narendra} and \cite{shorten} for a proof of 1) and 2).\hfill$\Box$

Now we do an inverse procedure to the Shorten-Narendra procedure, that is we have a $(n-1)\times (n-1)$ matrix and obtain a $n\times n$ matrix. Consider a $(n-1)\times (n-1)$ matrix $D_{n-1}$ defined as 
$$D_{n-1}=\begin{pmatrix}
d_{11} & d_{12} & \cdots & d_{1,n-1} \\
d_{21} & d_{22} & \cdots & d_{2,n-1} \\
\vdots & \vdots & \ddots & \vdots \\
d_{n-1,1} & d_{n-1,2} & \cdots & d_{n-1,n-1}
\end{pmatrix}$$
Then the matrix $A_n$ such that $A_n=
\begin{pmatrix}
A_{n-1} & d_{n-1} \\ c_{n-1}^T & d_{n-1}
\end{pmatrix}$ and $D_{n-1} = A_{n-1}- \frac{d_{n-1}c_{n-1}^T}{d_{n-1}}$ can be written as 
$$ A_n = \begin{pmatrix}
d_{11}-\frac{h_1 k_1}{k} &  \cdots & d_{1,n-1}-\frac{h_1 k_{n-1}}{k} & h_1 \\
d_{21}-\frac{h_2 k_1}{k} &  \cdots & d_{2,n-1}-\frac{h_2 k_{n-1}}{k} & h_2 \\ 
\vdots & \ddots & \vdots & \vdots \\
d_{n-1,1}-\frac{h_{n-1} k_1}{k} & \cdots & d_{n-1,n-1}-\frac{h_{n-1} k_{n-1}}{k} & h_{n-1} \\
k_1 & \cdots & k_{n-1} & -k
\end{pmatrix}. $$ 

Using above results and ideas presented in \cite{aguirre1} and \cite{aguirre2}, now we will prove that $\Sym\mathcal{H}_{n\times n}$ is diffeomorphic to $\Sym\mathcal{H}_{(n-1)\times(n-1)}\times\mathbb{R}^n$, where $\Sym\mathcal{H}_{n\times n}$ is the set of the Hurwitz symmetric matrices, actually that diffeomorphism will be an isomorphism of smooth fiber bundles, where the projection over $\Sym\mathcal{H}_{(n-1)\times(n-1)}$ is given by $A_n \longmapsto D_{n-1}$. Also we show the usefulness of these results for obtaining families of Hurwitz symmetric matrices or Hurwitz Metzler matrices and this let us do an analysis of robustness for important models in applications as the Insulin Model.

Another recent work related to ours, is presented in \cite{Briat:2017} where is addressed the sign-stability problem and it is proved the sign-stability of a Metzler sign-matrix can be assessed from the Hurwitz stability of a single particular matrix. Also, in Theorem 12.6 of \cite{Battacharyya} similar result about stability is stablished for a family of matrices in an interval matrix. On another hand, it follows from Theorem \ref{t1} (in next section) the existence of a neighborhood, for each $n\times n$-Hurwitz symmetric matrix, in which every matrix in such neighborhood is also a Hurwitz symmetric matrix, wich implies a more general robust stability without any reference to sign-matrix or special neighborhood. 

A different approach for the study of Hurwitz stability is given in \cite{SWS:2017} for linear continuous systems, through Schur complements and Lyapunov functions.

\section{Hurwitz symmetric matrices} \label{families-matrices}
As it was denoted in the introduction, let $\Sym\mathcal{H}_{n\times n}$ be the set of the Hurwitz symmetric matrices.

\begin{theorem}\label{t1} The space $\Sym\mathcal{H}_{n\times n}$ is diffeomorphic to the product manifold $\Sym\mathcal{H}_{(n-1)\times(n-1)}\times\mathbb{R}^n$.
\end{theorem}

\noindent\textbf{Proof.} First we define $$\varphi:\Sym\mathcal{H}_{n\times n}\rightarrow \Sym\mathcal{H}_{(n-1)\times(n-1)}\times\mathbb{R}^n.$$

If $A\in \Sym\mathcal{H}_{n\times n}$ with
\begin{equation}
A=\left(\begin{array}{cccc}a_{11}&a_{12}&\cdots&a_{1n}\\a_{21}&a_{22}&\cdots&a_{2n}\\\vdots&\vdots&\ddots&\vdots\\a_{n1}&a_{n2}&\cdots&a_{nn}\end{array}\right),\label{A}
\end{equation}
then $\varphi(A)$ is given by,
\begin{eqnarray} \label{A2} \left(\left(\begin{array}{ccc}a_{11}-\frac{a_{1n}a_{n1}}{a_{nn}}&\cdots&a_{1,n-1}-\frac{a_{1n}a_{n,n-1}}{a_{nn}}\\a_{21}-\frac{a_{2n}a_{n1}}{a_{nn}}&\cdots&a_{2,n-1}-\frac{a_{2n}a_{n,n-1}}{a_{nn}}\\ \vdots&\ddots&\vdots \\ a_{n-1,1}-\frac{a_{n-1,1}a_{n1}}{a_{nn}}&\cdots&a_{n-1,n-1}-\frac{a_{n-1,n}a_{n,n-1}}{a_{nn}}\end{array}\right),  a_{n1},a_{n2},\ldots,\ln(-a_{nn})\right).
\end{eqnarray}
Note that the matrix defined by \eqref{A} is Hurwitz symmetric if and only if, the matrix in the first coordinate in \eqref{A2} is Hurwitz symmetric matrix. Consequently, both of matrices have all negative entries on the diagonal.

Let $\varphi^{-1}:\Sym\mathcal{H}_{(n-1)\times(n-1)}\times\mathbb{R}^n\longrightarrow \Sym\mathcal{H}_{n\times n}$ be the inverse map of $\varphi$. If
$$ B=\left(\left(\begin{array}{cccc}b_{11}&b_{12}&\cdots&b_{1,n-1}\\b_{21}&b_{22}&\cdots&a_{2,n-1}\\\vdots&\vdots&\ddots&\vdots\\b_{n1}&b_{n2}&\cdots&b_{n-1,n-1}\end{array}\right),k_1,k_2,\ldots,k_{n-1},k_n\right), $$ with $B\in \Sym\mathcal{H}_{(n-1)\times(n-1)}\times\mathbb{R}^n$ then
$$ \varphi^{-1}(B)=\left(\begin{array}{cccc}b_{11}-\frac{k_1^2}{e^{k_n}}&\cdots&b_{1,n-1}-\frac{k_1k_{n-1}}{e^{k_n}}&k_1\\b_{21}-\frac{k_1k_2}{e^{k_n}}&\cdots&b_{2,n-1}-\frac{k_2k_{n-1}}{e^{k_n}}&k_2\\\vdots&\ddots&\vdots&\vdots\\b_{n-1,1}-\frac{k_1k_{n-1}}{e^{k_n}}&\cdots&b_{n-1,n-1}-\frac{k_{n-1}^2}{e^{k_n}}&k_{n-1}\\k_1&\cdots&k_{n-1}&-e^{k_n}\end{array}\right). $$ 

Is clear that $\varphi$ and $\varphi^{-1}$ are differentiable. Therefore, $\varphi$ is a diffeomorphism and we get the result we wanted to prove. $\hfill\Box$

\section{Application: Families of matrices}
When we are working with physical phenomena and we have to study a Hurwitz symmetric matrix $A$, it is more suitable to think that $A$ is an element of a family of matrices, since uncertainties in its elements must be considered. Then, we should be verify that all of the elements of a family of matrices are Hurwitz symmetric matrices. The results in this paper let us study complex families of matrices in $\Sym\mathcal{H}_{n\times n}$ which were constructed with simplest families in $\Sym\mathcal{H}_{(n-1)\times(n-1)}$. By example, if we have a ball $B_r$ with radius $r$, that is
\begin{equation*}
B_r = \left\{ \left. \left(
\begin{array}{cccc}
b_{11} & b_{12} & \cdots & b_{1,n-1} \\
b_{21} & b_{22} & \cdots & b_{2,n-1} \\
\vdots & \vdots & \ddots & \vdots    \\
b_{n1} & b_{n2} & \cdots & b_{n-1,n-1}
\end{array} \right) \right| \sum_{1\leq i,j \leq n-1} b_{ij}^2<r^2 \right\}
\end{equation*}
then, if $F:\Sym\mathcal{H}_{n\times n}\rightarrow \Sym\mathcal{H}_{(n-1)\times(n-1)}$, 

\begin{equation*}
F^{-1}(B_r) =   \left\{ \left. \left(
\begin{array}{cccc}
b_{11}-\frac{k_1^2}{e^{k_n}} & \cdots&b_{1,n-1}-\frac{k_1k_{n-1}}{e^{k_n}} & k_1 \\
b_{21}-\frac{k_1k_2}{e^{k_n}}& \cdots&b_{2,n-1}-\frac{k_2k_{n-1}}{e^{k_n}} & k_2 \\
\vdots & \ddots & \vdots & \vdots \\
b_{n-1,1}-\frac{k_1k_{n-1}}{e^{k_n}} & \cdots & b_{n-1,n-1}-\frac{k_{n-1}^2}{e^{k_n}}&k_{n-1} \\ 
k_1 & \cdots & k_{n-1} & -e^{k_n}
\end{array}
\right) \right| \sum_{1\leq i,j \leq n-1} b_{ij}^2<r^2; k_1,k_2,\ldots,k_n \in \mathbb{R} \right\}.
\end{equation*}

\noindent\textbf{Remark.} Here we use balls in the set of $\mathcal{M}_{(n-1)\times(n-1)}$ to show that our approach let us do an analysis of robustness for the stability of complex families of matrices in $\mathcal{M}_{n\times n}$.

\section{Metzler Hurwitz matrices}
\begin{definition}
The matrix $A$ is a strictly Metzler matrix, if its off diagonal elements are positive.
\end{definition}

We denote by $\Metz^+\mathcal{H}_{n\times n}$ the set of the $n\times n$ strictly Metzler Hurwitz matrices.

Consider the $(n-1)\times (n-1)$ Metzler Hurwitz matrix
\begin{equation} \label{B-matrix}
B= \left(\begin{array}{cccc}b_{11}&b_{12}&\cdots&b_{1,n-1}\\b_{21}&b_{22}&\cdots&b_{2,n-1}\\\vdots&\vdots&\ddots&\vdots\\b_{n1}&b_{n2}&\cdots&b_{n-1,n-1}\end{array}\right) 
\end{equation}
Now we look for all of $n\times n$ Metzler Hurwitz matrices $A= \left(\begin{array}{cc}A_{n-1} & b_{n-1} \\ c_{n-1}^T & d_{n-1} \end{array}\right)$ such that $B= A_{n-1} - \frac{b_{n-1} c_{n-1}^T }{d_{n-1}}$.

Then we can see that $A$ have to be written as
\begin{equation} \hspace*{-4mm}  \label{A-matrix}
A= \left(\begin{array}{cccc}
b_{11}-\frac{h_1 k_1}{e^{k_n}} & \cdots & b_{1,n-1}-\frac{h_1 k_{n-1}}{e^{k_n}} & h_1 \\ 
b_{21}-\frac{h_2 k_1}{e^{k_n}} & \cdots & b_{2,n-1}-\frac{h_2 k_{n-1}}{e^{k_n}} & h_2 \\ 
\vdots & \ddots & \vdots & \vdots \\ 
b_{n-1,1}-\frac{h_{n_1}k_1}{e^{k_n}} & \cdots & b_{n-1,n-1}-\frac{h_{n-1} k_{n-1}}{e^{k_n}} & k_{n-1} \\
k_{1} & \cdots & k_{n-1} & -e^{k_n} 
\end{array}\right)
\end{equation}
where $h_1, h_2, \ldots, h_{n-1}, k_1, k_2, \ldots, k_{n-1}, k_n \in \mathbb{R}$.

Note that $A$ is not a Metzler matrix for every $h_1, h_2, \ldots, h_{n-1}, k_1, k_2, \ldots, k_{n-1}, k_n \in \mathbb{R}$. We have to consider some restrictions on the $h_i$'s and $k_i$'s in order to get that $A$ is a Metzler matrix. This implies that the set of the $n\times n$ Metzler Hurwitz matrices is not a product manifold of the set of $(n-1)\times (n-1)$ Metzler Hurwitz matrix and the euclidean space. However, we can take advantages of our approach (as in Theorem \ref{t1}) for considering the \textit{fiber} of $(n-1)\times(n-1)$ Metzler Hurwitz matrices $B$ and then to get calculate families of $n\times n$ Metzler Hurwitz matrices and in this way we could study systems that consider uncertainties or perturbations. The result is the following.

\begin{theorem}
If $B$ defined in (\ref{B-matrix}) is a Metzler Hurwitz matrix, $k_i \geq 0$ for all $i=1, \ldots, n$, $h_j \geq 0$ for all $j=1, \ldots, n$ and $b_{ij}-\frac{h_i k_j}{e_{k_n}} \geq 0$ for $i\neq j$, then the matrix $A$ defined in (\ref{A-matrix}) is a Metzler Hurwitz matrix.
\end{theorem}

\noindent\textbf{Proof.} It follows of the above discussion and the results of Narendra and Shorten \cite{narendra} $\hfill\Box$

We illustrate this result with the following application.

\subsection{The Insulin Model}
Several mathematical models have been proposed for dynamic \textit{glucose-insulin therapy} in type 1 diabetes, however the mathematical model of Sorensen is one of the most widely accepted for its completeness in the representing of the metabolism of glucose, from a compartmental approach, see \cite{Sorensen}. The model is divided into three subsystems: glucose, insulin and rates of metabolic glucagon. The glucose subsystem is a nonlinear system of ordinary differential equations in eight dimensions, while the insulin subsystem is a linear system in seven dimensions. Both systems are coupled by nonlinear subsystem of rates of metabolic glucagon. The Insulin subsystem is $$ \dot x=Ax, $$ where $$A = \begin{pmatrix}
-\frac{173}{100} & \frac{173}{100} & 0 & 0 & 0 & 0 & 0 \vspace*{1mm}\\
\frac{227}{500} & -\frac{3151}{1000} & 0 & \frac{909}{1000} & \frac{727}{1000} & \frac{53}{50} & 0 \vspace*{1mm}\\
0& \frac{153}{200} & -\frac{153}{200} & 0 & 0 & 0 & 0 \vspace*{1mm}\\
0& \frac{47}{500} & \frac{189}{500} & -\frac{789}{1000} & 0 & 0 & 0 \vspace*{1mm}\\
0& \frac{1411}{1000} & 0 & 0 & -\frac{367}{200} & 0 & 0 \vspace*{1mm}\\
0& \frac{709}{500} & 0 & 0 & 0 & -\frac{937}{500} & \frac{91}{200} \vspace*{1mm}\\
0& 0 & 0 & 0 & 0 & \frac{1}{20} & -\frac{111}{1000}
\end{pmatrix}
$$ and $x=(x_1,\ x_2\cdots,\  x_7)^T$, where the entries of this vector are the insulin concentrations in the compartment of brain, heart, stomach, liver, kidney, vascular space and interstitium of the peripheral tissue; respectively, see \cite{quiroz}. The nominal parameters used here are those validated in \cite{Sorensen}.
If we renamed $A$ as $A_7$, then,
$ A_7 = \begin{pmatrix} A_6 & b_6 \\ c_6^T & d_6 \end{pmatrix}$, we can calculate $B_6 = A_6-\frac{b_6 c_6^T}{d_6}$,
$$B_6 = \begin{pmatrix}
-\frac{173}{100} & \frac{173}{100} & 0 & 0 & 0 & 0 \vspace*{1mm}\\
\frac{227}{500} & -\frac{3151}{1000} & 0 & \frac{909}{1000} & \frac{727}{1000} & \frac{53}{50} \vspace*{1mm}\\
0& \frac{153}{200} & -\frac{153}{200} & 0 & 0 & 0 \vspace*{1mm}\\
0& \frac{47}{500} & \frac{189}{500} & -\frac{789}{1000} & 0 & 0 \vspace*{1mm}\\
0& \frac{1411}{1000} & 0 & 0 & -\frac{367}{200} & 0 \vspace*{1mm}\\
0& \frac{709}{500} & 0 & 0 & 0 & -\frac{937}{500} + \frac{91}{444}
\end{pmatrix}
$$

We have that the fiber of $B_6$ is the familie of $7\times 7$-matrices given by

\begin{equation}\label{fiberB6}
\begin{pmatrix}
-\frac{173}{100}-\frac{h_1 k_1}{e^{k_7}} & \frac{173}{100}-\frac{h_1 k_2}{e^{k_7}} & -\frac{h_1 k_3}{e^{k_7}}& -\frac{h_1 k_4}{e^{k_7}}& -\frac{h_1 k_5}{e^{k_7}}& -\frac{h_1 k_6}{e^{k_7}} & h_1 \vspace*{1mm}\\
\frac{227}{500}-\frac{h_2 k_1}{e^{k_7}} & -\frac{3151}{1000}-\frac{h_2 k_2}{e^{k_7}} & -\frac{h_2 k_3}{e^{k_7}}& \frac{909}{1000}-\frac{h_2 k_4}{e^{k_7}} & \frac{727}{1000}-\frac{h_2 k_5}{e^{k_7}} & \frac{53}{50}-\frac{h_2 k_6}{e^{k_7}} & h_2 \vspace*{1mm}\\
-\frac{h_3 k_1}{e^{k_7}}& \frac{153}{200}-\frac{h_3 k_2}{e^{k_7}} & -\frac{153}{200}-\frac{h_3 k_3}{e^{k_7}} & -\frac{h_3 k_4}{e^{k_7}}& -\frac{h_3 k_5}{e^{k_7}}& -\frac{h_3 k_6}{e^{k_7}} & h_3 \vspace*{1mm}\\
-\frac{h_4 k_1}{e^{k_7}}& \frac{47}{500}-\frac{h_4 k_2}{e^{k_7}} & \frac{189}{500}-\frac{h_4 k_3}{e^{k_7}} & -\frac{789}{1000}-\frac{h_4 k_4}{e^{k_7}} & -\frac{h_4 k_5}{e^{k_7}}& -\frac{h_4 k_6}{e^{k_7}} & h_4 \vspace*{1mm}\\
-\frac{h_5 k_1}{e^{k_7}}& \frac{1411}{1000}-\frac{h_5 k_2}{e^{k_7}} & -\frac{h_5 k_3}{e^{k_7}}& -\frac{h_5 k_4}{e^{k_7}}& -\frac{367}{200}-\frac{h_5 k_5}{e^{k_7}} & -\frac{h_5 k_6}{e^{k_7}} & h_5 \vspace*{1mm}\\
-\frac{h_6 k_1}{e^{k_7}}& \frac{709}{500}-\frac{h_6 k_2}{e^{k_7}} & -\frac{h_6 k_3}{e^{k_7}}& -\frac{h_6 k_4}{e^{k_7}}& -\frac{h_6 k_5}{e^{k_7}}& -\frac{937}{500} + \frac{91}{444}-\frac{h_6 k_6}{e^{k_7}} & h_6 \vspace*{1mm}\\
k_1 & k_2 & k_3 & k_4 & k_5 & k_6 & e^{-k_7}  
\end{pmatrix}
\end{equation}
where $k_i \in \mathbb{R}$, $i=1,\ldots, 7$, $h_j \in \mathbb{R}$ and $j=1,\ldots, 6$.

Note that the matrices defined in (\ref{fiberB6}) are not Metzler for every value of the $k_i$'s and $h_j$'s. If we give the conditions $k_1=0$, $k_2 \geq 0$, $k_3=0$, $k_4=0$, $k_5=0$, $k_6 \geq 0$, $k_7 \in \mathbb{R}$, $h_1=0$, $h_2 \geq 0$, $h_3=0$, $h_4=0$, $h_5=0$ and $h_6 \geq 0$, then we get the family of matrices
\begin{equation}\label{particular-family}
\begin{pmatrix}
-\frac{173}{100} &\hspace{-2mm} \frac{173}{100} &\hspace{-2mm} 0 &\hspace{-2mm} 0 &\hspace{-2mm} 0 &\hspace{-2mm} 0 &\hspace{-2mm} 0 \vspace*{1mm}\\
\frac{227}{500} &\hspace{-1mm}\hspace{-2mm} -\frac{3151}{1000}-\frac{h_2 k_2}{e^{k_7}} &\hspace{-2mm} 0 &\hspace{-2mm} \frac{909}{1000} & \frac{727}{1000} & \frac{53}{50}-\frac{h_2 k_6}{e^{k_7}} & h_2 \vspace*{1mm}\\
0&\hspace{-2mm} \frac{153}{200} &\hspace{-2mm} -\frac{153}{200} &\hspace{-2mm} 0 &\hspace{-2mm} 0 &\hspace{-2mm} 0 &\hspace{-2mm} 0 \vspace*{1mm}\\
0&\hspace{-2mm} \frac{47}{500} &\hspace{-2mm} \frac{189}{500} &\hspace{-2mm} -\frac{789}{1000} &\hspace{-2mm} 0 &\hspace{-2mm} 0 &\hspace{-2mm} 0 \vspace*{1mm}\\
0&\hspace{-2mm} \frac{1411}{1000} &\hspace{-2mm} 0 &\hspace{-2mm} 0 &\hspace{-2mm} -\frac{367}{200} &\hspace{-2mm} 0 &\hspace{-2mm} 0 \vspace*{1mm}\\
0&\hspace{-2mm} \frac{709}{500}-\frac{h_6 k_2}{e^{k_7}} &\hspace{-2mm} 0 &\hspace{-2mm} 0 &\hspace{-2mm} 0 &\hspace{-2mm} -\frac{23158}{13875}-\frac{h_6 k_6}{e^{k_7}} &\hspace{-2mm} h_6 \vspace*{1mm}\\
0& k_2 & 0 & 0 & 0 & k_6 & -e^{k_7}
\end{pmatrix}
\end{equation}

If we take $k_2$, $k_6$, $h_2$ and $h_6$ nonnegative, with $k_7 \in \mathbb{R}$ such that $\frac{709}{500}-\frac{h_6 k_2}{e^{k_7}} \geq 0$ and $\frac{53}{50}-\frac{h_2 k_6}{e^{k_7}}\geq0$ then all of the elements of the familie (\ref{particular-family}) are Metzler matrices.

Since $B_6$ is a Hurwitz matrix, then all of the elements of (\ref{particular-family}) are Hurwitz matrices.

Note that $A_7$ is in the familie (\ref{particular-family}): take $h_2=k_2=0$, $k_6=\frac{1}{20}$, $h_6=\frac{91}{200}$, $k_7=\ln(\frac{111}{1000})$. 

The importance of (\ref{particular-family}) is that now we have an Insulin Model that consider uncertainties or perturbations. Therefore, now we can answer questions like the following; if nominal values of the matrix $A$ components correspond to the real system (the human body of the insulin dependent person), \textit{any statement about the stability of the system can be altered if we change one of the inputs?}, since in this section we have shown a robustness result for a family of systems of linear differential equations. Furthermore, considering the inoculation of insulin in the above linear model, this becomes a positive linear systems of the form $\dot x=Ax+b$ with $A$ Metzler matrix and $b\in\mathbb{R}^n_+$, it is known (cf. \cite{BNS:89}) that if $A$ is also a Hurwitz matrix is a necessary and sufficient condition to there exists a no trivial positive fixed point $\overline x\in\mathbb{R}^n_+$, such that, $\overline x=-A^{-1}b$, this means that if the dynamic is disturbed, all solutions of the linear system (including the fixed point) does not leave the positive cone.


\section{Conclusions}
In this paper, we present an approach of differential topology for studying robust stability of the set of Hurwitz-Metzler matrices and Hurwitz symmetric matrices. In a specific problem, the analysis of robust stability can be done with standard technics of linear algebra, but our approach let us do the analysis for families of complex systems, as the Remark at the end of Section \ref{families-matrices} indicated.

\section*{Acknowledgements}
Jes\'us F. Espinoza acknowledges the financial support of CONACyT and of the Universidad de Sonora.

\end{document}